\documentclass[10pt,twocolumn,twoside,journal]{IEEEtran} 

\usepackage{psfrag}
\usepackage{graphicx}
\usepackage{cite}
\usepackage{amsmath} 
\usepackage{amsthm} 
\usepackage{amssymb}  
\usepackage{amsfonts}
\usepackage{mathrsfs}
\usepackage{mathtools}
\usepackage[dvipsnames]{xcolor}
\usepackage{enumerate}
\usepackage{circuitikz}
\usepackage{pgfplots}
\usepackage{booktabs}
\usepackage{bbm}
\usepackage[bbgreekl]{mathbbol}
\usepackage{colortbl}

\usepackage{etoolbox}

\newcommand{\beq}{\begin{equation}}
\newcommand{\eeq}{\end{equation}}
\newcommand{\bseq}{\begin{subequations}}
\newcommand{\eseq}{\end{subequations}}
\newcommand{\beqn}{\begin{eqnarray}}
\newcommand{\eeqn}{\end{eqnarray}}
\newcommand{\ba}{\begin{array}}
\newcommand{\ea}{\end{array}}
\newcommand{\bct}{\begin{center}}
\newcommand{\ect}{\end{center}}
\newcommand{\btab}{\begin{tabular}}
\newcommand{\etab}{\end{tabular}}
\newcommand{\btmz}{\begin{itemize}}
\newcommand{\etmz}{\end{itemize}}
\newcommand{\benum}{\begin{enumerate}}
\newcommand{\eenum}{\end{enumerate}}

\newcommand{\R}{{\mathbb R}}

\newcommand{\norm}[1]{\| #1 \|}                 

\newcommand{\matbegin}{
        \left[
}
\newcommand{\matend}{
        \right]
}

\newcommand{\tbo}[2]{
  \matbegin \begin{array}{c}
       #1 \\ #2
       \end{array} \matend }

\newcommand{\tbt}[4]{
  \matbegin \begin{array}{cc}
       #1 & #2 \\ #3 & #4
       \end{array} \matend }

\newcommand{\be}{\begin{equation}}
\newcommand{\ee}{\end{equation}}

\newcommand{\cplxs}{ C\kern -.35em \rule{0.03 em}{.7 ex}~   }

\def\complex{\hbox{C\kern -.45em \rule{0.03 em}{1.5 ex}}~}

\newcommand{\bi}{\begin{itemize}}
\newcommand{\ei}{\end{itemize}}

\newtheorem{assumption}{Assumption}

\DeclareSymbolFontAlphabet{\mathbb}{AMSb}
\DeclareSymbolFontAlphabet{\mathbbl}{bbold}

\DeclareMathOperator*{\argmin}{argmin}

\providecommand{\norm}[1]{\left\|#1\right\|}

\newcommand{\mc}{\mathcal}

\newcommand{\one}{\mathbf{1}}

\newcommand{\bbR}{\mathbb{R}}

\newtheorem{theorem}{Theorem}
\newtheorem{lemma}{Lemma}
\newtheorem{proposition}{Proposition}
\newtheorem{remark}{Remark}
\newtheorem{definition}{Definition}
\newtheorem{corollary}{Corollary}

\graphicspath{{Figures/} }

\pgfplotsset{compat=1.14}

\title{Online Optimization as a Feedback Controller:\\ Stability and Tracking}
\author{Marcello Colombino, Emiliano Dall'Anese and Andrey Bernstein\thanks{The authors are with the National Renewable Energy Laboratory (NREL), Golden, CO, USA. Email: \{name.lastname\}@nrel.gov. The work was supported by the DOE Office of Electricity, Advanced Grid Modeling Program.}} 

\begin{document}

\maketitle

\begin{abstract}
This paper develops and analyzes feedback-based online optimization methods to regulate the output of a linear time-invariant (LTI) dynamical system to the optimal solution of a \emph{time-varying} convex  optimization problem. The design of the algorithm is based on continuous-time primal-dual dynamics, properly modified to incorporate feedback from the  LTI dynamical system, applied to a proximal augmented Lagrangian function. The resultant closed-loop algorithm tracks  the solution of the time-varying optimization problem without requiring knowledge of (time-varying) disturbances in the dynamical system.  The analysis leverages integral quadratic constraints to provide linear matrix inequality (LMI) conditions that guarantee global exponential stability and  bounded tracking error. Analytical results show that, under a sufficient time-scale separation between the dynamics of the LTI dynamical system and the algorithm, the LMI conditions can be always satisfied. The paper further  proposes a modified algorithm that can track an approximate solution trajectory of the constrained optimization problem under less restrictive assumptions. As an illustrative example, the proposed algorithms are showcased for power transmission systems, to compress the time scales between secondary and tertiary control, and allow to simultaneously power re-balancing  and tracking of DC optimal power flow points. 
\end{abstract}

\section{Introduction}

The problem of regulating the output of a dynamical system to optimal solutions of an optimization problem is common to many engineering systems~\cite{Jokic09,Elia-CDC11,Brunner-CDC12}. The traditional modus operandi in e.g., chemical processes, power systems and energy networks, as well as  many other cyber-physical systems involves a time-scale separation between optimization and regulation/control tasks. In this setting, optimization tasks produce reference signals at a slow(er) time-scale, based on well-defined operational costs, engineering constraints, and algebraic representations of the system input-output map. Reference signals are then fed to reference-tracking controllers, which drive in real time the operation of the system to the prescribed reference points while coping with disturbances. A typical example of time-scale separation can be found in power transmission systems: fast local controllers act at the (sub)-second level to maintain power balance, while at a much slower time-scale (i.e., every 15 minutes) the system operator solves a given optimization problem based on forecasts of loads and renewable generations, a grid model, and security constraints to provide set-points for the controllable assets~\cite{Wollenbergbook}. 

A different approach to controlling a system to an economical and feasible state is to formulate the control problem directly as a constrained finite-horizon optimization problem and capitalize on model predictive control methods~\cite{garcia1989model}. While the second approach has clear advantages such as constraint satisfaction at each time-step, both methodologies rely on system models and accurate  forecasts of the disturbances affecting the system.

Motivated by e.g., flow congestion control in communication systems~\cite{low1999optimization} and lately by real-time control of power systems~\cite{bacher1988real,commelec1,Bolognani_feedback_15,anese2016optimal,hauswirth2016projected,tang2017real,bernstein2018online,Mazzi2018} and machine learning (e.g.,~\cite{shalev2012online,TianyiChen2017_Bandits} and pertinent references therein), \emph{feedback-based optimization} algorithms exploits measurements form the system to steer the dynamical systems to the solution of an optimization problem with limited knowledge of the input-output map and no information about the disturbance. In particular, feedback-based algorithms were developed in~\cite{anese2016optimal,bernstein2018online} to track solutions of time-varying problems. In these works, however, the algorithms are designed based on  an algebraic map, and underlying system dynamics are neglected. 

This paper proposes online optimization methods to regulate the output of a dynamical system to the solution of a \emph{time-varying convex optimization problem}, while taking into account the dynamics of the underlying dynamical system.  In particular, linear time-invariant (LTI) dynamical system are considered. The design principles involve continuous-time primal-dual gradient methods, applied to a proximal augmented Lagrangian function, as a \emph{feedback controller}. Stability and tracking of the proposed online feedback optimization algorithm leverages robust control methods based on integral quadratic constraints (IQCs)~\cite{megretski1997system}. Proximal augmented Lagrangian method for (open-loop) continuous saddle-flows and IQCs were recently pioneered by~\cite{dhingra2016proximal}. IQCs have proven to be very effective to study also the convergence of discrete-time first-order optimization algorithms; see e.g.,~\cite{lessard2016analysis,fazlyab2017dynamical}.  Similar techniques can be found in e.g.,~\cite{nelson2017integral}, where the authors propose a co-design of a state estimator and an online optimizer for unconstrained optimization problems. Lastly, this paper differs from~\cite{Rahili2015} by considering time-varying constrained problems and not restricting the design to integrator dynamics. Overall, the contributions of this paper are the following: 

\noindent \emph{(i)} We provide linear matrix inequality (LMI)  conditions that guarantee global exponential tracking bounds for the system output with respect of the solution of a time-varying constrained optimization problem; 

\noindent \emph{(ii)} We show that, under sufficient time-scale separation between the LTI dynamical system and the online algorithm, the LMIs are always feasible and therefore the online optimizer is stable and drives the system to the optimal solution; and, 

\noindent \emph{(iii)} We propose and analyze an \emph{approximate} online optimization algorithm that provably converges to an approximate optimal trajectory under less restrictive assumptions.

To conclude the paper, an example of application  for the proposed  algorithms in the power systems context is provided. In particular, the paper considers a DC optimal power flow (OPF) problem an considers the transmission system swing dynamics. It is shown that the online optimization algorithm provides a means to compress the time scales between secondary frequency control and constrained economic dispatch. The proposed methodology subsumes the frameworks of~\cite{li2016connecting,dorfler2016breaking}, and extends then to incorporate more general time-varying convex optimization settings. The proposed algorithms allow cost- and constraint-aware system rebalancing without any information of  uncontrollable loads and uncertain renewable generation.

\noindent \emph{Notation}. The set of real numbers is denoted by $\R$. $\R_+$ denotes the set of nonnegative reals.
Given a matrix $A$, $A^\top$ denotes its transpose. We use $\bar{\sigma}(A)$ to indicate the largest singular value of $A$, $\underline{\sigma}(A)$ its smallest singular value, $\kappa(A)$ its condition number and $\mathrm{null}(A)$ its null space. We write $A\succeq0$  $(A\succ0)$ to denote that $A$ is symmetric and positive semidefinite (definite). {For column vectors $x\in\R^n$ and $y\in\R^m$ we use $(x,y) = [x^\top, y^\top]^\top \in \R^{n+m}$ to denote a stacked vector.} Furthermore, $I$ denotes the identity matrix, $\mathbf 1$ is the vector of all ones, $\norm{x}$ denotes the Euclidean norm, and, given a matrix $P\succ0$ we define $\norm{x}_P=\sqrt{x^\top P x}$.





\section{The online optimization framework}\label{sec.framework}

This paper focuses on control methods to regulate the output of a dynamical system to the optimal solution of a time-varying optimization problem. This section outlines the problem formulation and the proposed solution method. 

\subsection{Problem setup}

Consider an LTI dynamical system of the form 
\begin{align}\label{sys}
\begin{split}
 \dot x &= Ax + B u + B_w w\\
y_1 &= C_1x + D_{1w}w\\
y_2 &= C_2x + D_{2w}w
\end{split}
\end{align}
where $A$ is Hurwitz stable; the vector $x$ represents the system's state; $u$ is the control input; $w$ is an exogenous disturbance; and,  $y$ is  the measurable output (the reason for separating $y$ into $y_1$ and $y_2$ will become apparent shortly). We define the transfer function of system~\eqref{sys} from $u$ to $y_i$ as $G_{iu}(s):=C_i(sI-A)^{-1}B$ and the transfer function from $w$ to $y_i$ as $G_{iw}(s):=C_i(sI-A)^{-1}B_w + D_{iw}$, for $i\in\{1,2\}$. Note that, for constant $u$ and $w$, $y_i$ is given by $y_i = G_{iu}(0)u + G_{iw}(0) w $. For convenience, let us denote the steady-state maps by
\[
\Pi_{iu} := G_{iu}(0),\quad \Pi_{iw} := G_{iw}(0),\quad i\in\{1,2\}.
\]

In many engineering applications, it is desirable to optimize over the steady-state of a dynamical system to ensure that the system settles to a safe and/or economically viable operating point. Consider a \emph{time-varying} constrained optimization problem over the steady-state of system~\eqref{sys} of the form:
\begin{align}\label{opt}
\begin{split}
\min_{u, y_1, y_2} 			&~f^t(u) + h^t(y_1) + g^t(y_2)\\
\text{s.t.}   &~y_1 = \Pi_{1u}u + \Pi_{1w} w(t) \\
		&~y_2 = \Pi_{2u}u + \Pi_{2w} w(t) 
\end{split}
\end{align}
where, for all $t\ge0$, $f^t$ is $m$-strongly convex with $L_f$-Lipschitz continuous gradient, $h^t$ is convex and has $L_h$-Lipschitz continuous gradient and $g^t$ is proper, convex and lower semi-continuos but not necessarily differentiable. The function  $f^t$ represents a time-varying cost function on the steady-state control effort, $h^t$ represents a cost on the steady-state system output e.g., tracking penalty or soft constraints and $g^t$ be the indicator function of a time-varying convex \emph{constraint set} or a non-differentiable regularizer. 
{Note that, in order to motivate problem~\eqref{opt} we inherently assume that the disturbance $w(t)$ and the functions $f^t, h^t$ and $g^t$ vary slower than the system~\eqref{sys}, so that the optimization over the steady state of \eqref{sys} is justified.}
In most engineering applications, problem~\eqref{opt} is solved off-line and the optimal solution $u^\star$ is fed to the systems~\eqref{sys} as a reference. This, however, requires knowledge or a forecast of the disturbance $w(t)$ and of the effect of the disturbance on the system's outputs $\Pi_w$. In this paper, we study a method to control the system~\eqref{sys} so that its output \emph{tracks} an optimal trajectory of problem~\eqref{opt} \emph{without} any knowledge of $w$ or $\Pi_w$.

\subsection{Proximal augmented Lagrangian algorithms}
In this section, we review a first-order algorithm~ developed in \cite{dhingra2016proximal} to solve \emph{time invariant} convex  problems of the form:
\begin{align}\label{opt.ti}
\begin{split}
\min_{u, y_1, y_2} 			&~f(u) + h(y_1) + g(y_2)\\
\text{s.t.}   &~y_1 = \Pi_{1u}u + \Pi_{1w} w \\
		&~y_2 = \Pi_{2u}u + \Pi_{2w} w 
\end{split}
\end{align}
Given $\mu>0$, we define the proximal operator of the function $g$ as
\[
\mathbf{prox}_{\mu g}(v):=\argmin_x g(x)+\frac{1}{2\mu}\|x-v\|^2,
\]
and the associated Moreau envelope of $g$ as
\[
M_{\mu g}(v):=g(\mathbf{prox}_{\mu g}(v))+\frac{1}{2\mu}\|\mathbf{prox}_{\mu g}(v)-v\|^2.
\] Note that if $g$ is proper and convex $M_{\mu g}$ is not only convex but also continuously differentiable and
\[
\nabla M_{g\mu} (v)=\frac{1}{\mu}(v-\mathbf{prox}_{\mu g}(v)). 
\] 
We refer the reader to~\cite{parboy13} for a comprehensive review of the properties of proximal operators. Next, we introduce the augmented Lagrangian for problem~\eqref{opt}
\begin{align}\label{eq.aug.lag}
\begin{split}
\mc L_{\mu}(u,\lambda):=& f(u) + h(\Pi_{1u}u+\Pi_{1w}w)\\
& + M_{\mu g}(\Pi_{2u}u+\Pi_{2w}w+\mu \lambda) - \frac{\mu}{2}\|\lambda\|^2
\end{split}
\end{align}
where $y_1$ is substituted in the cost function and the constraints on $y_2$ is dualized. The saddle points of $\mc L_{\mu}(u,\lambda)$ are precisely the optimal primal-dual pairs of problem~\eqref{opt.ti}~\cite{dhingra2016proximal}. A version of the augmented Lagrangian~\eqref{eq.aug.lag} was  introduced in~\cite{rockafellar1976augmented} when $g$ represents a convex constraint set and has recently been adopted in~\cite{qu2018exponential} for algorithmic developments. In~\cite{dhingra2016proximal}, the authors introduce the proximal augmented Lagrangian in its general form~\eqref{eq.aug.lag}, and propose the saddle-flow algorithm\footnote{In order to recover exactly the setup of~\cite{dhingra2016proximal}, one can define a function $g^w(v):=g(v+\Pi_{2w}w)$ and define the augmented Lagrangian with respect to this new function. Given the properties of the proximal operator~\cite[Section 2.2]{parboy13} the algorithm remains unchanged.}
\begin{align*}
\begin{split}
 \dot u = - \nabla_u \mc L_{\mu}(u,\lambda),\quad \dot \lambda = \nabla_\lambda \mc L_{\mu}(u,\lambda),
 \end{split}
\end{align*}
which, for Problem~\eqref{opt.ti} becomes
\begin{align}\label{offline.opt}
\begin{split}
 \dot u = &  -\nabla f(u) -\Pi_{1u}^\top\nabla h(\Pi_{1u}u+\Pi_{1w}w)\\
 & -\Pi_{2u}^\top \nabla M_{g\mu} (\Pi_{2u}u+\Pi_{2w}w+\mu\lambda) \\
 \dot \lambda   = & \,\mu\bigg( \nabla M_{g\mu} (\Pi_{2u}u+\Pi_{2w}w+\mu\lambda) -\lambda \bigg).
 \end{split}
\end{align}
and, given the strong convexity assumption, they prove that the algorithm~\eqref{offline.opt} globally converges to the \emph{unique} primal-dual pair of problem~\eqref{opt.ti} and, under further technical assumptions has an exponential convergence rate, which is in contrast to standard saddle-flow algorithms for constrained problems~\cite{cherukuri2016role} without the proximal augmented Lagrangian.

In our setting, computing $\nabla h(\Pi_{1u}u+\Pi_{1w}w)$ and  $\nabla M_{g\mu}(\Pi_{2u}u+\Pi_{2w}w+\mu\lambda)$ requires information about the disturbance $w$ and the disturbance model $\Pi_w$, which is impractical in many engineering applications. Furthermore, the optimization problem is often time-varying and therefore it is desirable to devise algorithms that can \emph{track} the optimal solution. This motivates the online implementation of the optimization algorithm~\eqref{offline.opt} in the next section.

\subsection{Online implementation}

In this section we consider an online implementation of the primal-dual saddle flow algorithm~\eqref{offline.opt} for the time-varying optimization problem~\eqref{opt}, in feedback with the dynamical system~\eqref{sys}. The aim of the feedback interconnection is achieving output tracking of the time-varying optimal trajectory of systems~\eqref{opt}. The main advantage is that the online version of the primal-dual method~\eqref{offline.opt} is implementable without any information on the disturbance $w$ or any knowledge on how $w$ affects the output. In fact, in the online implementation, the effect of $w$ is ``observed" through the measurement of the system's output $y$.
The online implementation of a time varying adaptation of the primal dual algorithm~\eqref{offline.opt} in feedback with the dynamical system~\eqref{sys} is given by
\begin{align}\label{online.opt}
\begin{split}
 \dot x &= Ax +   B u +   B_w w\\
 \dot u &= -\nabla f^t(u)- \Pi_{1u}^\top\nabla h^t(y_1) -\Pi_{2u}^\top \nabla M_{g^t\mu} (y_2+\mu\lambda) \\
 \dot \lambda &= \mu\bigg( \nabla M_{g^t\mu} (y_2+\mu\lambda) -\lambda \bigg)\\
  y_1 &= C_1x + D_{1w}w\\
 y_2 &= C_2x + D_{2w}w,
\end{split}
\end{align}
where the time dependence of $x,u,y_1,y_2$ and $w$ was dropped for ease of notation. Let us formalize the assumptions needed for the theoretical results.
\begin{assumption}\label{ass.conv}
For all $t\ge0$ The function $f^t$ is $m$-strongly convex with $L_f$-Lipschitz continuous gradient, $h^t$ is convex and has $L_h$-Lipschitz continuous gradient and $g^t$ is proper, convex and lower semi-continuous.
\end{assumption}
Next, we characterize the equilibria of the closed-loop system~\eqref{online.opt}.
\begin{definition}\label{def.stat} Given a time-varying dynamical system $\dot x = F^t(x,w)$ and a time-varying disturbance $w(t)$, a time-varying stationary point is a set valued map defined by $\Phi(t):=\{\bar x\,|\, 0 = F^t(\bar x,w(t))\}$.
\end{definition}
\begin{proposition}\label{propo.optimality}
Given a time-varying disturbance $w(t)$, if problem~\eqref{opt} is feasible at time $t$, the set $\Phi(t)$ of  time-varying stationary points of system~\eqref{online.opt} according to Definition~\ref{def.stat} is a singleton of the form $(-A^{-1}Bu^\star(t), u^\star(t), \lambda^\star(t))$, where $(u^\star(t),\lambda^\star(t))$ is the unique primal-dual optimal pair at time $t$ for problem~\eqref{opt}.
\end{proposition}
\begin{proof}
According to Definition~\ref{def.stat}, $(\bar x,\bar u, \bar \lambda)$ is a time-varying stationary point of~\eqref{online.opt} if
\begin{subequations}
\begin{align}\label{online.opt.eq}
 \bar x &= -A^{-1}B\bar u + -A^{-1}B_w w(t)\\
 \begin{split}
 0 &= -\nabla f^t(\bar u -\Pi_{1u}^\top \nabla h^t(\Pi_{1u}\bar u+\Pi_{1w}w(t)\\
 & -\Pi_{2u}^\top \nabla M_{g^t\mu} (\Pi_{2u}\bar u +\Pi_{2w}w(t) +\mu\bar  \lambda \label{eq.u}  \end{split} \\
 0 &= \mu\bigg( \nabla M_{g^t\mu} (\Pi_{2u}\bar u + \Pi_ww(t) +\mu\bar \lambda -\bar \lambda \bigg). \label{eq.lambda}
\end{align}
\end{subequations}
It is simple to verify that the stationarity conditions of~\eqref{eq.u}-\eqref{eq.lambda} are equivalent to the first order optimality conditions of~\eqref{opt}~\cite{dhingra2016proximal}.
\end{proof}
Proposition~\ref{propo.optimality}, establishes that, for each $t\ge 0$, the (time-varying) stationary point of the closed loop system~\eqref{online.opt} corresponds to the (time-varying) primal-dual optimal pair of~\eqref{opt}. Thus, if the online scheme~\eqref{online.opt} is stable, the proposed control strategy has the ability to ``steer" the output of the dynamical system~\eqref{sys} towards the (time-varying) optimal solution of problem~\eqref{opt}.
Note that the proposition refers to the time varying stationary point $\bar x(t)$ of the linear system~\eqref{sys} and \emph{not} the state evolution $x(t)$.
Finally, two technical assumptions are presupposed.
\begin{assumption}\label{ass.feas}
The time varying problem~\eqref{opt} is feasible for all $t\ge0$ and the function $(u^\star,\lambda^\star):[0,\infty)\to\R^m\times \R^p$, which maps time to the optimal primal-dual trajectory of~\eqref{opt}, is measurable.
\end{assumption}
\begin{assumption}\label{ass.cont}
The RHS of~\eqref{online.opt} is measurable and continuous for almost all $t\ge0$.
\end{assumption}
These assumptions hold naturally for most reasonable time varying functions $f^t,h^t,g^t$, e.g. piece-wise constant. Providing conditions that guarantee the satisfaction of the assumptions is beyond the scope of this paper.
Assumptions~\ref{ass.conv}-\ref{ass.cont} allow $f^t,g^t,h^t$ and $w(t)$ to vary discontinuously but still guarantee the existence of a Caratheodory solution to system~\eqref{online.opt}~\cite[Proposition S1]{cortes2008discontinuous}. A Caratheodory solution of a dynamical system of the form
$
\dot x = f(x,t),
$
is an absolutely continuous map $x:[0,\infty)\to\R^n$ that solves the (Lebesgue) integral version of the differential equation, i.e.,
\[
x(t) = x(0) + \int_0^tf(x(\tau),\tau)\mathrm{d}\tau.
\]
or, equivalently, solves the differential equation for almost all $t$. We refer the reader to~\cite{cortes2008discontinuous} for a comprehensive review on dynamical systems with discontinuous right hand side.

The online optimization scheme~\eqref{online.opt} is de-facto a feedback control law for system~\eqref{sys} designed to \emph{track} the optimal solution of problem~\eqref{opt}. In the following we analyze in detail the stability properties and tracking performance  of~\eqref{online.opt}.

\section{Stability and performance analysis}\label{sec.stability}

This section contains the main results of the paper. We analyze the online optimization scheme~\eqref{online.opt} and we provide a tractable stability test to verify that, under the proposed control law, the closed-loop system~\eqref{online.opt} is stable and the output of system~\eqref{sys} tracks the optimal solution of problem~\eqref{opt}. 

\subsection{Tracking optimal trajectories for time-varying disturbances}

\begin{figure}[t]
\begin{center}
\input{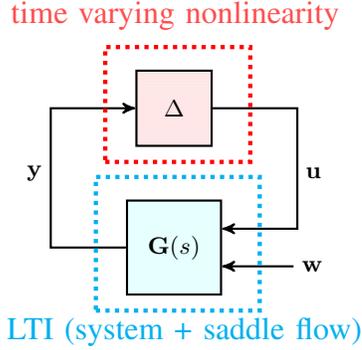}%
\noindent%
%
\begin{tikzpicture}[scale=1, auto, >=stealth']

    \small


     \node[clpblock, minimum height=1.25cm, minimum width=1.25cm] (N)  {$\mathbf G(s)$};

     \node[deltablock, minimum height=1cm, minimum width=1cm, anchor=south]
                      (delta) at ($(N.north) + (0,0.7cm)$) {$\Delta$};

     \def\ioOffset{0.25cm}
     \coordinate (Nv) at ($(N.east) + (0,\ioOffset)$);
     \coordinate (Nu) at ($(N.east) - (0,\ioOffset)$);
     \coordinate (Nz) at ($(N.west) + (0,0)$);
     \coordinate (Ny) at ($(N.west) - (0,\ioOffset)$);

    \node (Numid) at ($(Nu) + (1.cm,0)$){};
    \node (Numid2) at ($(Numid) + (0.2cm,0)$){$\mathbf w$};
    \node (Nymid) at ($(Ny)-(1.cm,0)$){};

    \draw [connector] (Nz) -- (Nz -| Nymid) |- (delta);
    \draw [connector]  (delta) -- (delta -| Numid) |- (Nv);

    \draw [connector]  (Numid2) --  (Nu);

    \path (N.north) -- node (midpoint) [midway]{} (delta.south) {};

    \node[left] at (Nymid |- midpoint) {$\mathbf y $};

    \node[right] at (Numid |- midpoint) {$\mathbf u$};

    \draw[red,ultra thick,dotted] ($(delta.north west)+(-0.4,0.3)$)  rectangle ($(delta.south east)+(0.5,-0.2)$);
    \node [red,label={[xshift=-0.8cm, yshift=0.3cm]{\large \textcolor{red!70}{time varying nonlinearity}}}] at ($(delta.north east)+(0.3,0.0)$) {};

     \draw[ProcessBlue,ultra thick,dotted] ($(N.north west)+(-0.4,0.3)$)  rectangle ($(N.south east)+(0.5,-0.2)$);
    \node [red,label={[xshift=-0.8cm, yshift=-0.9cm]{\large \textcolor{ProcessBlue}{LTI (system + saddle flow)}}}] at ($(N.south east)+(0.3,0.0)$) {};

\end{tikzpicture} 
\caption{The online optimization algorithm can be written as a feedback interconnection of a LTI system $\mathbf G(s)$ and a time varying nonlinearity $\Delta$. $\mathbf G(s)$ contains both the physical system and part of the optimization routine while $\Delta$ contains information on the gradients of the time-varying objective function and constants. }
\label{fig.gdelta}
\end{center}
\end{figure}

As a first step, we re-write system~\eqref{online.opt} as the feedback interconnection of a LTI system with a ``benign" time varying nonlinearity as illustrated in Figure~\ref{fig.gdelta}. The linear system $\mathbf G(s)$ is given by
 \begin{align}\label{online.opt.iqc}
\begin{split}
\dot {\mathbf{z}} & = \mathbf A \mathbf{z} + \mathbf{B}\mathbf{u} + \mathbf{B_w}\mathbf{w}\\
 \mathbf{y} & = \mathbf C \mathbf{z} + \mathbf{ D_w} \mathbf w\\
\mathbf u & = \Delta (\mathbf y)
\end{split}
\end{align}
where $\mathbf z :=(x,u,\lambda)$, $\mathbf w :=w$ the matrices $\mathbf {A},\mathbf {B},\mathbf {C}\text { and }\mathbf {D}$ are defined as 
 \begin{align}\label{eq.matrices}
 \begin{split}
\mathbf A := 
\begin{bmatrix}  A\!\!&\!\!  B \!\!& \!\!0\\
0 \!\!&\!\! -mI \!\!&\!\! 0\\
0 \!\!&\!\! 0  \!\!&\!\! -\mu I
 \end{bmatrix}\!, 
\mathbf{B}:= \begin{bmatrix}
0 \!\!& \!\!0 \!\!& \!\!0\\
-I \!\!& \!\!-\Pi_{1u}^\top \!\!& \!\!-\frac{1}{\mu}\Pi_{2u}^\top\\
0 \!\!& \!\!0 \!\!&  \!\!I
 \end{bmatrix}\!,\\
 \mathbf{B_w}:= \begin{bmatrix}
B_w \\
0\\
0
 \end{bmatrix}\!,
\mathbf {C} := \begin{bmatrix}
0 & I & 0\\
C_1 & 0 & 0\\
C_2 & 0 & \mu I\\
 \end{bmatrix}\!,
 \mathbf {D_w} := \begin{bmatrix}
0\\ 
D_{1w}\\
D_{2w}
 \end{bmatrix}\!,
 \end{split}
\end{align}
and the nonlinear function $\Delta(\cdot)$ is defined as
 \begin{align*}
 {\Delta(\mathbf y)} := 
\begin{bmatrix} 
\nabla f^t(\mathbf y_1) - m \mathbf y_1 \\
\nabla h^t(\mathbf y_2)\\
\mu \nabla M_{g^t\mu} (\mathbf{y}_3)
 \end{bmatrix}
\end{align*}
\begin{definition}\label{def.opt}
 Given any $\mathbf w(t)$, $\mathbf {z^\star_w}(t)$ is an \underline{optimal trajectory} if 
 \begin{align}\label{online.opt.ss}
\begin{split}
0 & = \mathbf A \mathbf {z^\star_w}(t) + \mathbf{B}\mathbf {u^\star_w}(t) + \mathbf{B_w}\mathbf{w}(t),
\end{split}
\end{align}
where
\begin{align*}
\begin{split}
 \mathbf{y_w^\star}(t) & = \mathbf C \mathbf {z^\star_w}(t) + \mathbf{ D_w} \mathbf w(t)\\
\mathbf{u^\star_w}(t) & = \Delta (\mathbf {y^\star_w}(t)). 
\end{split}
\end{align*}
\end{definition}
In the remainder of the document, in the interest of readability, we will drop the $t$ dependency of the optimal trajectory unless explicitly needed. From Proposition~\ref{propo.optimality}, we conclude that $\mathbf {z^\star_w}$ is of the form $\mathbf {z^\star_w} = (-A^{-1}Bu^\star,u^\star,\lambda^\star)$, where $(u^\star,\lambda^\star)$ is a primal-dual optimal pair for problem~\eqref{opt}. Since $\mathbf w$ is time-varying, $\mathbf {z^\star_w}$ is also \emph{time-varying} and represents the optimal \emph{trajectory} of problem~\eqref{opt}. In the following we give an LMI condition on the matrices $\mathbf{A,B,C}$ that guarantee that the online optimization scheme~\eqref{online.opt} is stable and tracks time-varying optimal trajectories of~\eqref{opt} exponentially fast. Given $\varphi_1,\varphi_2\ge0$ and $\hat L_f:= L_f - m$, we define
\[
\Xi_\varphi:=\begin{bmatrix} 
0 & 0 & 0 & \varphi_1\hat L_f I & 0& 0\\
0 & 0 & 0 & 0 & \varphi_2L_h I & 0\\
0 & 0 & 0 & 0 & 0& I\\
\varphi_1 \hat L_f I & 0 & 0& -2\       I &0 & 0\\
0 & \varphi_2 L_h I  & 0 & 0 &  -2\varphi_2I  & 0\\
0 & 0 & I & 0 & 0 & -2I\\
 \end{bmatrix}
\]

\begin{definition}\label{def.lmi}
Given the matrices $\mathbf{(A,B,C,D)}$ and $(\mathbf P,\Xi_\varphi)$, we say that $\mathbf{(A,B,C,D)}\in \mathbb{}{L}(\mathbf{P},\Xi_\varphi)$ if
\begin{align}\label{LMI}
\begin{bmatrix} 
\mathbf{A^\top P + P A}   & \mathbf{PB} \\
 \mathbf{B^\top P} & 0 
 \end{bmatrix}
+ 
\begin{bmatrix} 
\mathbf{C^\top}   & 0 \\
0 & I 
 \end{bmatrix}
\Xi_\varphi
 \begin{bmatrix} 
\mathbf{C}   & 0 \\
0 & I 
 \end{bmatrix}\preceq 0
\end{align}
\end{definition}

\begin{theorem}\label{thm.kyp}
Let $\mathbf w$ be a time-varying disturbance and  $\mathbf{A}_{\rho}:= \mathbf A+\rho I$, under Assumptions~\ref{ass.conv}-\ref{ass.cont}, if there exist $\rho>0, \varphi_1,\varphi_2\ge0$ and matrix $\mathbf P\succ 0$ such that $\mathbf{(A_\rho,B,C,D)}\in \mathbb{}{L}(\mathbf{P},\Xi_\varphi)$ according to Definition~\ref{def.lmi}, then
\begin{align*}
\begin{split}
\norm{\mathbf{z}(t)-\mathbf{z^\star_w}(t)} \le& \sqrt{\kappa(\mathbf P)}\norm{\mathbf{z}(0)-\mathbf{z^\star_w}(0)}\mathrm e^{-\rho t} \\
&+ \frac{\kappa(\mathbf P)}{2} \int_0^t e^{-\rho (t-\tau)} \norm{\mathbf{\dot z^\star_w(\tau)}}\mathrm d \tau,
\end{split}
\end{align*}
where $\kappa(\mathbf P)$ is the condition number of $\mathbf P$, $\mathbf z(t)$ is the trajectory of system~\eqref{online.opt} and $\mathbf {z^\star_w}(t)$ is the optimal trajectory According to Definition~\ref{def.opt} and $\dot{\mathbf{z}}_{\mathbf{w}}^\star(t)$ is its distributional derivative, i.e., $\int_0^t\dot{\mathbf{z}}_{\mathbf{w}}^\star(\tau)\mathrm{d}\tau = {\mathbf{z_w}}^\star(t)$.
\end{theorem}
\begin{proof}
Since, for all $t\ge0$, $\nabla f^t(v) - m  v$ is the $\hat L_f-$Lipschitz continuous gradient of the convex function $f^t(v)-\frac{m}{2}\|v\|^2$ with $\hat L_f= L-m$, $\nabla h^t(v)$ is the $L_h-$Lipschitz continuous gradient of the convex function $h^t(v)$, and $\mu \nabla M_{g^t\mu} (v)$ is the scaled 1-Lipschitz continuous gradient of the convex Moreau envelope $M_{g^t\mu} (v)$, the nonlinear function $\Delta(\cdot)$ satisfies the following point-wise-in-time IQC defined by $\Xi_\varphi$, that is, for any arbitrary $\mathbf {\hat y}(t), \mathbf {\hat u}(t)$ such that $\mathbf {\hat u}(t)=\Delta(\mathbf {\hat y}(t))$  and $\varphi_1,\varphi_2 \ge 0$~\cite[Lemma 6]{lessard2016analysis}
 \begin{align}\label{eq.iqc}
   \begin{split}
\begin{bmatrix} 
\mathbf y - \mathbf {\hat y}(t)\\
\mathbf u- \mathbf {\hat u}(t)
 \end{bmatrix}^\top
\Xi_\varphi
\begin{bmatrix} 
\mathbf y - \mathbf {\hat y}(t)\\
\mathbf u- \mathbf {\hat u}(t)
 \end{bmatrix}&\ge 0, \\
  \forall \,\,(\mathbf y,\mathbf u) \text{ s.t. } \mathbf u=\Delta(\mathbf y),\quad \forall~ t&\ge0.
  \end{split}
\end{align}
Consider the now time-varying Lyapunov function $V(\mathbf{z})= \norm{\mathbf{ z}-\mathbf{ z}^\star_{\mathbf w}}_{\mathbf P}^2 $. We pre- and post-multiply the matrices in~\eqref{LMI} by $(\mathbf{ z}-\mathbf{ z}^\star_{\mathbf w}, \mathbf{u}-\mathbf{ u}^\star_{\mathbf w})$, noting that
 \begin{align*}
 \ba{rcl}
\dot {\mathbf{z}} & \!= \!&\mathbf A (\mathbf{z}-\mathbf{z^\star_w}) + \mathbf{B}(\mathbf{u}-\mathbf{u^\star_w}) \\
 \mathbf{y} - \mathbf{y^\star_w}& \!= \! & \mathbf C(\mathbf{z}-\mathbf{z^\star_w}) ,
\ea
\end{align*}
we conclude that the derivative of $V$ along the trajectories of~\eqref{online.opt.iqc} satisfies
\begin{align}\label{eq.derivative.w}
\begin{split}
\dot V &\le- 2\rho V - 
\begin{bmatrix}
\mathbf{y} - \mathbf{y^\star_w}\\
\mathbf{u}-\mathbf{u^\star_w}
\end{bmatrix}^\top
\Xi_\varphi
\begin{bmatrix}
\mathbf{y} - \mathbf{y^\star_w}\\
\mathbf{u}-\mathbf{u^\star_w}
\end{bmatrix} + (\mathbf{z}-\mathbf{z^\star_w})\mathbf P\mathbf{\dot z^\star_w}\\
&\le -2\rho V + \bar\sigma(\mathbf P)\norm{\mathbf{z}-\mathbf{z^\star_w}} \norm{\mathbf{\dot z^\star_w}},
\end{split}
\end{align}
where the last term derives from the fact that the optimal trajectory is time varying; notice that this is the key difference of our analysis compared to the standard analysis of non-time-varying systems.
Now consider the function $W = \sqrt V = \norm{\mathbf{z}-\mathbf{z^\star_w}}_{\mathbf P}$. Since $\dot W = \frac{\dot V}{2\sqrt V}$ we conclude that
\begin{align}\label{eq.derivative.w2}
\begin{split}
\dot W &\le-\rho W + \frac{\bar\sigma(\mathbf P)}{2\sqrt{\underline\sigma(\mathbf P)}} \norm{\mathbf{\dot z^\star_w}}.
\end{split}
\end{align}
By the comparison lemma~\cite{K02}, the following norm inequality is satisfied
\begin{align*}
\begin{split}
\norm{\mathbf{z}(t)-\mathbf{z^\star_w}(t)}_{\mathbf P} \le& \norm{\mathbf{z}(0)-\mathbf{z^\star_w}(0)}_{\mathbf P}\mathrm e^{-\rho t} \\
&+ \frac{\bar\sigma(\mathbf P)}{2\sqrt{\underline\sigma(\mathbf P)}} \int_0^t e^{-\rho (t-\tau)} \norm{\mathbf{\dot z^\star_w}}\mathrm d \tau.
\end{split}
\end{align*}
and therefore
\begin{align*}
\begin{split}
\norm{\mathbf{z}(t)-\mathbf{z^\star_w}(t)} \le& \sqrt{\kappa(\mathbf P)}\norm{\mathbf{z}(0)-\mathbf{z^\star_w}(0)}\mathrm e^{-\rho t} + \\
&\frac{\kappa(\mathbf P)}{2} \int_0^t e^{-\rho (t-\tau)} \norm{\mathbf{\dot z^\star_w}}\mathrm d \tau,
\end{split}
\end{align*}
which completes the proof.
\end{proof}
Theorem~\ref{thm.kyp} establishes an LMI condition to verify exponential tracking of the optimal trajectory for system~\eqref{online.opt} and can be interpreted in as an Input to State Stability (ISS) result~\cite{sontag2008input}, where the input to the system is the rate of change of the optimizer. With the following two corollaries we recover some results in the literature for the discrete-time static case (online optimization in feedback to a static map), which has recently received increasing attention in the context of power systems~\cite{bernstein2018online,tang2017real,hauswirth2016projected}.
\begin{corollary}\label{coro.exp}
If the disturbance $\mathbf w$ is \emph{constant in time}, then, if $\mathbf{(A_\rho,B,C,D)}\in \mathbb{}{L}(\mathbf{P},\Xi_\varphi)$, system~\eqref{online.opt} exponentially converges to the unique (and constant in time) optimal solution $\mathbf{z^\star_w}$ satisfying~\eqref{online.opt.ss}, i.e.,
\begin{align*}
\begin{split}
\norm{\mathbf{z}(t)-\mathbf{z^\star_w}} &\le \sqrt{\kappa(\mathbf P)}\norm{\mathbf{z}(0)-\mathbf{z^\star_w}}\mathrm e^{-\rho t}.
\end{split}
\end{align*}
\end{corollary}
Note that, given the IQC formulation of the problem, Corollary~\ref{coro.exp} would follow from~\cite[Theorem 3]{hu2016exponential}. 
Finally, the next corollary covers  the case in which $\norm{\mathbf{\dot z^\star_w}}$ is bounded and provides an asymptotic bound for the tracking error.
\begin{corollary}~\label{coro.limsup}
Suppose $\norm{\mathbf{\dot z^\star_w}}\le\sigma$ for all $t>0$, then, if $\mathbf{(A_\rho,B,C,D)}\in \mathbb{}{L}(\mathbf{P},\Xi_\varphi)$, 
\begin{align*}
\begin{split}
\limsup_{t\to\infty} \norm{\mathbf{z}(t)-\mathbf{z^\star_w}(t)} &\le \frac{\sigma\kappa(\mathbf P)}{2\rho}.
\end{split}
\end{align*}
\end{corollary}
Corollary~\ref{coro.limsup} is the dynamical counterpart of the main result in~\cite{bernstein2018online} for online optimization in feedback with a static map.


\subsection{Time-scale separation and feasibility of the LMI condition}

Theorem~\ref{thm.kyp} guarantees that, if the LMI test~\eqref{LMI} is satisfied, the online feedback optimization scheme~\eqref{online.opt} is exponentially stable and, in case of a time varying disturbance, tracks the optimal trajectory of the time-varying optimization problem~\eqref{opt}. In the following, we prove that, under an assumption on the rank of $\Pi_{2u}$, if the system is fast enough compared to the optimization algorithm, there always exist $\rho>0$, $\varphi_1$, $\varphi_2$, and $\mathbf{P} \succ 0$ such that $\mathbf{(A_\rho,B,C,D)}\in \mathbb{}{L}(\mathbf{P},\Xi_\varphi)$ according to Definition~\ref{def.lmi}. 
To this end, we consider a ``fast" version of system~\eqref{sys} given by
\begin{align}\label{sys.fast}
\begin{split}
\epsilon\, \dot x &= Ax + Bu + B_w w\\
y_1 &= C_1x + D_{1w}w\\
y_2 &= C_2x + D_{2w}w
\end{split}
\end{align}

and define 
$$
\mathbf A_{\epsilon} := 
\begin{bmatrix} \frac{1}{\epsilon} A& \frac{1}{\epsilon} B & 0\\
0 & -mI & 0\\
0 & 0  & -\mu I
 \end{bmatrix}
$$
and 
$\mathbf{A}_{\epsilon\rho}:= \mathbf A_{\epsilon}+\rho I$. 
Let us define 
\[
\bar \Pi_{1u}:= \begin{bmatrix} 
I\\
\Pi_{1u}
 \end{bmatrix}
\]
and $m_\rho := m -\rho$, $\mu_\rho := \mu -\rho$. 
Note that the steady-state of system~\eqref{sys.fast} is the same as that of system~\eqref{sys}. This means that the optimal trajectory $\mathbf z^\star$ of the online optimization algorithm applied to the faster system remains unchanged. The following technical lemma is needed for the proof of the main result of the section.
\begin{lemma}\label{lem.limit}
Let $\mathbf G_{\epsilon\rho}(s):=\mathbf C\left(sI-\mathbf A_{\epsilon\rho} \right)^{-1}\mathbf B$, with $\mathbf{B,C}$ defined in~\eqref{eq.matrices}, then
$
\lim_{\epsilon\to0}\mathbf G_{\epsilon\rho}(s) = \mathbf H_{\rho}(s),
$
where 
\begin{align*}
\mathbf H_{\rho}(s)& := \begin{bmatrix}
\bar \Pi_{1u} & 0\\
\Pi_{2u}& \mu I
\end{bmatrix}
\begin{bmatrix}g_{m\rho}(s) I\!\! &\!\! 0 \\
0\!\! &\!\! g_{\mu\rho}(s) I
\end{bmatrix}
\begin{bmatrix}
 -\bar \Pi_{1u}^\top & \!\!-\frac{1}{\mu}\Pi_{2u}^\top\\
  0 &   \!\!I
\end{bmatrix}
\end{align*}
and $g_{m\rho}(s):=(s-m_\rho)^{-1}$ and $g_{\mu\rho}(s):=(s-\mu_\rho)^{-1}$.
\end{lemma}
\begin{proof}
We note that $\left(sI-\mathbf A_{\epsilon\rho} \right)^{-1}$ is given by
\begin{align*}
\begin{bmatrix} (sI-\frac{1}{\epsilon} A +\rho I)^{-1} \!\!& \!\! \frac{1}{\epsilon}(sI-\frac{1}{\epsilon} A +\rho I)^{-1} B g_{m\rho}(s) \!\!& \!\!0\\
0 \!\!&\!\! g_{m\rho}(s)I\!\!&\!\! 0\\
0 \!\!&\!\! 0 \!\! &\!\! g_{\mu\rho}(s)I
 \end{bmatrix}
\end{align*}
and $\lim_{\epsilon\to0}\frac{1}{\epsilon}C_i (sI-\frac{1}{\epsilon} A+\rho I)^{-1} B = -C_iA^{-1}B = \Pi_{iu}$. The result follows from simple algebraic manipulations. 
\end{proof}
We are now ready to present the main result of this section
\begin{theorem}\label{thm.time.scale.separation}
If $\Pi_{2u}\Pi_{2u}^\top\succ 0$ and $\mu\ge\max\{\hat L_f,L_h\}$, under Assumptions~\ref{ass.conv}-\ref{ass.cont} there exist an $\epsilon>0$ and $\rho>0,\varphi_1,\varphi_2\ge0$ and a matrix $\mathbf P\succ 0$ such that  $\mathbf{(A_{\epsilon\rho},B,C,D)}\in \mathbb{}{L}(\mathbf{P},\Xi_\varphi)$.
\end{theorem}
\begin{proof}
For fixed $\epsilon,\rho,\varphi_1,\varphi_2$, using the KYP Lemma~\cite{megretski1997system} we conclude that condition~\eqref{LMI} holds for $
\mathbf A_{\epsilon\rho}$ if and only if
\begin{align}\label{FDI}
\begin{bmatrix} 
\mathbf G_{\epsilon\rho}(j\omega)\\
I
 \end{bmatrix}^\star
 \Xi_\varphi
 \begin{bmatrix} 
\mathbf G_{\epsilon\rho}(j\omega)\\
I
 \end{bmatrix}
\preceq 0,\quad \forall \, \omega\in \R
\end{align}
Using Lemma~\ref{lem.limit} we know that
\begin{align*}
\lim_{\epsilon\to 0}
\begin{bmatrix} 
\mathbf G_{\epsilon\rho}(j\omega)\\
I
 \end{bmatrix}^\star
 \Xi_\varphi
 \begin{bmatrix} 
\mathbf G_{\epsilon\rho}(j\omega)\\
I
 \end{bmatrix}
= 
\begin{bmatrix} 
\mathbf H_{\rho}(j\omega)\\
I
 \end{bmatrix}^\star
 \Xi_\varphi
 \begin{bmatrix} 
\mathbf H_{\rho}(j\omega)\\
I
 \end{bmatrix}
\end{align*}

Similarly to~\cite[Theorem 4]{dhingra2016proximal}, we can prove that, under the assumptions of the theorem, there exist $\rho^\star>0$, $\varphi_1,\varphi_2$ and $\mu$, such that
\begin{align}\label{FDI2}
\begin{bmatrix} 
\mathbf H_{\rho}(j\omega)\\
I
 \end{bmatrix}^\star
 \Xi_\varphi
 \begin{bmatrix} 
\mathbf H_{\rho}(j\omega)\\
I
 \end{bmatrix}
\preceq 0,\quad \forall \, \omega\in \R.
\end{align}
Let us define $\bar L\ge\max\{\hat L_f ,L_h\}$ and the new multiplier 
\[
\bar{\Xi}_\varphi:=\begin{bmatrix} 
0 & 0 & 0 & \varphi_1\bar L I & 0& 0\\
0 & 0 & 0 & 0 & \varphi_2\bar L I & 0\\
0 & 0 & 0 & 0 & 0& I\\
\varphi_1 \bar L I & 0 & 0& -2\varphi_1 I &0 & 0\\
0 & \varphi_2 \bar L I  & 0 & 0 &  -2\varphi_2I  & 0\\
0 & 0 & I & 0 & 0 & -2I\\
 \end{bmatrix}
\]
Note that if
\begin{align}\label{FDI3}
\begin{bmatrix} 
\mathbf H_{\rho}(j\omega)\\
I
 \end{bmatrix}^\star
\bar \Xi_\varphi
 \begin{bmatrix} 
\mathbf H_{\rho}(j\omega)\\
I
 \end{bmatrix}
\preceq 0,\quad \forall \, \omega\in \R,
\end{align}
then~\eqref{FDI2} is satisfied. This is because using $ \bar\Xi_\varphi$ guarantees robust stability with respect to a \emph{larger} class of functions $\Delta(\cdot)$ as, by definition, $\bar L-$Lipschitz function are also $\hat L_f-$ and $L_h-$ Lipschitz. Let us choose $\varphi_1=\varphi_2=1$, $\mu = \bar L$, 
with simple manipulations we infer that~\eqref{FDI3} holds for the particular choice of parameters if and only if, for all $\omega\in \R$
 \begin{align}\label{FDI4}
	\tbt{\!\!\!\dfrac{\mu m_\rho}{m_\rho^2 + \omega^2}  \, \bar \Pi_{1u}\bar \Pi_{1u}^\top + I\!\!\!\!\!}
	{\!\!\!\!\!\dfrac{m_\rho}{m_\rho^2 + \omega^2} \,\bar \Pi_{1u} \Pi_{2u}^\top\!\!\!}
	{\!\!\!*\!\!\!\!\!}
	{\!\!\!\!\!\dfrac{m_\rho/\mu}{m_\rho^2 + \omega^2} \,\Pi_{2u}\Pi_{2u}^\top + \dfrac{\omega^2 - \rho \mu_\rho}{\mu_\rho^2 + \omega^2} \, I\!\!\!}
	\!\succeq\!
	0.
\end{align}
Consider now the case in which $\rho=0$. Clearly, since 
$\Pi_{2u}\Pi_{2u}^\top\succ0$, the bottom-right block in~\eqref{FDI4} is positive definite. Therefore, we consider its Schur complement, which can be bounded by
\[
S_{22_{|\rho=0}}\succeq I+ \bar \Pi_{1u}\left(\dfrac{\mu m}{m^2 + \omega^2}-\tilde\Pi_{2u} \dfrac{\mu m}{ (m^2 + \omega^2)}\right)\bar \Pi_{1u}^\top 
\]
where $\tilde\Pi_{2u}: = \Pi_{2u}(\Pi_{2u}\Pi_{2u}^\top)^{-1}\Pi_{2u}^\top$ is a projection matrix. Since $0\preceq\tilde\Pi_{2u}\preceq I$, we conclude that  
$
S_{22_{|\rho=0}}\succeq  I \succ0,
$
which implies that~\eqref{FDI4} holds \emph{strictly} for $\rho=0$. By continuity, there exist $\rho^\star>0$ such that~\eqref{FDI4} holds strictly for all $\rho<\rho^\star$. With a similar continuity argument, there exist $\epsilon^\star>0$ such that, for all $\epsilon<\epsilon^\star$,~\eqref{FDI} is satisfied with $0\le\rho<\rho^\star$, which completes the proof.
\end{proof}
\begin{remark}\label{rem.time.scale}
Note that, if the system
$
\epsilon \dot x = f(x,z),\, \dot z = g(x,z)
$
is exponentially stable with parameter $\rho$ then
$
 \dot x = f(x,z),\,\dot z = \epsilon g(x,z)
$
is exponentially stable with parameter $\epsilon\rho$.
\end{remark}
In view of Remark~\ref{rem.time.scale}, provided that the rank condition $\Pi_{2u}\Pi_{2u}^\top\succ0$ is satisfied, Theorem~\ref{thm.time.scale.separation} provides a \emph{constructive} tool for designing a stabilizing online optimizer by artificially slowing down the optimization algorithm. Also note that Theorem~\ref{thm.time.scale.separation} guarantees that, if there is enough time-scale separation, the LMI condition~\eqref{LMI}, which is at the base of Theorem~\ref{thm.kyp}, is feasible; it does not however say that time-scale separation is necessary. In many practical examples,~\eqref{LMI} is still satisfied even though there is no apparent time-scale separation between the system and the optimizer. 
The rank condition $\Pi_{2u}\Pi_{2u}^\top\succ0$, which is necessary for Theorem~\ref{thm.time.scale.separation}, implies that, in problem~\eqref{opt}, we can constrain at most as many outputs as there are inputs using \emph{hard constraints} through the indicator function $g$. This is the reason of the particular formulation of~\eqref{opt}, in which we can use \emph{soft constraints} through $h^t$. It is observed in practice that the rank constraint is \emph{necessary} for the LMI condition in~\eqref{LMI} to be satisfied with $\rho>0$. Proving the role of this assumption in the exponential stability of the online optimization~\eqref{online.opt} is left as future work. In the next section we proposed an approximate algorithm to overcome this assumption.

\section{Approximate online optimization}\label{sec.approx}

In order to overcome the need of the rank condition $\Pi_{2u}\Pi_{2u}^\top\succ0$, we propose the following variation of the algorithm~\eqref{online.opt} with regularizing parameter $\gamma$.

\begin{align}\label{online.opt.gamma}
\begin{split}
 \dot x &= Ax +   Bu +   B_w w\\
 \dot u &= -\nabla f^t(u)- \Pi_{1u}^\top\nabla h^t(y_1) -\Pi_{2u}^\top \nabla M_{g^t\mu} (y_2+\mu\lambda) \\
 \dot \lambda &= \mu\bigg( \nabla M_{g^t\mu} (y_2+\mu\lambda) -\lambda \bigg) - \gamma Q\lambda\\
  y_1 &= C_1x + D_{1w}w\\
 y_2 &= C_2x + D_{2w}w
\end{split}
\end{align}
Where $Q\succeq0$ is selected such that $\mathrm{null}(\Pi_{2u}\Pi_{2u}^\top)\cap\mathrm{null}(Q)=\{0\}$. 
The regularized algorithm~\eqref{online.opt.gamma} can be used as an alternative to soft constraints in the case in which  $\Pi_{2u}\Pi_{2u}^\top\not\succ0$. For example, suppose $u\in\R^m$, $y_2\in\R^p$ with $p>m$. If $g(y_2)$ defines the indicator function of a convex time-varying constraint set $\mc Y^t= \mc Y_1^t\times\dots\times\mc Y_m^t$, then $\mathbf{prox}_{\mu g}(v) = \mathrm{Proj}_{\mc Y^t}(v)$ and the stationarity condition on the multiplier $\lambda$ for the un-regularized algorithm~\eqref{online.opt} imply
\begin{align*}
  0 & = \mu\bigg( \nabla M_{g^t\mu} (y_2^\star+\mu\lambda^\star) -\lambda^\star \bigg)\\
 \iff 0 & = \mu\bigg( \nabla M_{g^t\mu} (y_2^\star+\mu\lambda) -\lambda^\star \bigg)\\
 \iff 0 & =   y_2^\star - \mathrm{Proj}_{\mc Y}(y_2^\star+\mu\lambda^\star) \implies y_2^\star  \in \mc Y.
\end{align*}
However, since $m>p$, we conclude that $\Pi_{2u}\Pi_{2u}^\top\not\succ0$. Therefore, we cannot use the control algorithm~\eqref{online.opt}. In that case we can choose $Q = \mathrm{diag}(0_p,I_m)$ and use the regularized algorithm~\eqref{online.opt.gamma}. The new stationarity conditions for the approximate algorithm read 
\begin{align*}
 & 0  = \mu\bigg( \nabla M_{g^t\mu} (y_2^{\gamma\star}+\mu\lambda^{\gamma\star}) -\lambda^{\gamma\star} \bigg) - \gamma Q\lambda^{\gamma\star} \\
 \iff &0  = \mu\bigg( \nabla M_{g^t\mu} (y_2^{\gamma\star}+\mu\lambda^{\gamma\star}) -\lambda \bigg)-\gamma Q\lambda^{\gamma\star}\\
 \iff& 0  =   y_2 - \mathrm{Proj}_{\mc Y^t}(y_2^{\gamma\star}+\mu\lambda^{\gamma\star}) \\
 \implies & y_2^{\gamma\star} + \gamma Q\lambda^{\gamma\star}  \in \mc Y^t.
\end{align*}The particular choice of $Q$ still guarantees that the first $p$ constraints are enforced strictly and the last $m-p$ are enforced approximately. The level of approximation can be adjusted with the parameter $\gamma$. This formulation is an alternative to soft constraints and can be beneficial in practical problems.
\subsection{Stability analysis of the approximate algorithm}
In this section we define the tracking properties of the approximate algorithm~\eqref{online.opt.gamma}. First, we define the \emph{approximate} optimal trajectory as
 \begin{align}\label{online.opt.ss.gamma}
\begin{split}
0 & = \mathbf {A^\gamma} \mathbf {z^{\star\gamma}_w} + \mathbf{B}\mathbf {u^{\star\gamma}_w} + \mathbf{B_w}\mathbf{w}\\
 \mathbf{y_w^{\star\gamma}} & = \mathbf C \mathbf {z^{\star\gamma}_w} + \mathbf{ D_w} \mathbf w\\
\mathbf{u^{\star\gamma}_w} & = \Delta (\mathbf {y^{\star\gamma}_w}),
\end{split}
\end{align}
where $\mathbf {A^\gamma}$ is defined as
\[
\mathbf {A^\gamma} := 
\begin{bmatrix}  A& B & 0\\
0 & -mI & 0\\
0 & 0  & -\mu I-\gamma Q
 \end{bmatrix}
\]
From Theorem~\ref{thm.kyp}, if there exists a $\rho>0$ such that  $\mathbf{(A_{\rho}^{\gamma\top},B,C,D)}\in {L}(\mathbf{P},\Xi_\varphi)$ for some $\mathbf P\succ 0$, $\varphi_1,\varphi_2\ge 0$,
 
\begin{align*}
\begin{split}
\norm{\mathbf{z}(t)-\mathbf{z^{\star\gamma}}(t)} \le& \sqrt{\kappa(\mathbf P)}\norm{\mathbf{z}(0)-\mathbf{z^{\star\gamma}_w}(0)}\mathrm e^{-\rho t} \\
&+ \frac{\kappa(\mathbf P)}{2} \int_0^t e^{-\rho (t-\tau)} \norm{\mathbf{\dot z^{\star\gamma}_w}}\mathrm d \tau.
\end{split}
\end{align*}
Next we show that, given enough time-scale separation, the partially regularized algorithm~\eqref{online.opt.gamma} can stabilize the system without the rank constraint $\Pi_{2u}\Pi_{2u}^\top\succ0$. To this end we define the matrices
\[
\mathbf {A_{\epsilon}^\gamma} := 
\begin{bmatrix} \frac{1}{\epsilon} A& \frac{1}{\epsilon} B & 0\\
0 & -mI & 0\\
0 & 0  & -\mu I-\gamma Q
 \end{bmatrix}
\]
and $\mathbf{A_{\epsilon\rho}^\gamma}:=\mathbf {A_{\epsilon}^\gamma}-\rho I.$
\begin{theorem}\label{thm.time.scale.separation.gamma}
If $\mathrm{null}(\Pi_{2u}\Pi_{2u}^\top)\cap\mathrm{null}(Q)=\{0\}$ and $\mu\ge\max\{\hat L_f,L_h\}$, under Assumptions~\ref{ass.conv} and~\ref{ass.feas}, for all $\gamma>0$, there exist an $\epsilon>0$ and $\rho>0,\varphi_1,\varphi_2\ge0$ and a matrix $\mathbf P\succ 0$ such that  $\mathbf{(A_{\epsilon\rho}^{\gamma\top},B,C,D)}\in {L}(\mathbf{P},\Xi_\varphi)$.
\end{theorem}
\begin{proof}
Let us define $\mathbf G^\gamma_{\epsilon\rho}(s):=\mathbf C(sI - \mathbf {A^\gamma_{\epsilon\rho}}^{-1}\mathbf{B}$, By the KYP Lemma, there exist an $\epsilon>0$ and $\rho>0,\varphi_1,\varphi_2\ge0$ and a matrix $\mathbf P\succ 0$ such that  $(\mathbf {A^\gamma_{\epsilon\rho},B,C,D})$ satisfy the LMI~\eqref{LMI} if and only if
\begin{align}\label{FDI.epsilon}
\begin{bmatrix} 
\mathbf G^\gamma_{\epsilon\rho}(j\omega)\\
I
 \end{bmatrix}^\star
 \Xi_\varphi
 \begin{bmatrix} 
\mathbf G^\gamma_{\epsilon\rho}(j\omega)\\
I
 \end{bmatrix}\preceq 0,
\end{align}
Using the same arguments as in the proof of Theorem~\eqref{thm.time.scale.separation}, we conclude that, for the particular choice if $\mu = \hat L$ $\varphi_1=\varphi_2=1$,
\begin{align*}
\lim_{\epsilon\to 0}
\begin{bmatrix} 
\mathbf G^\gamma_{\epsilon\rho}(j\omega)\\
I
 \end{bmatrix}^\star
 \Xi_\varphi
 \begin{bmatrix} 
\mathbf G^\gamma_{\epsilon\rho}(j\omega)\\
I
 \end{bmatrix}\preceq 0,
\end{align*}
if and only if
\begin{align}\label{FDI5}
	\tbt{\!\!\!\dfrac{\mu m_\rho}{m_\rho^2 + \omega^2}  \, \bar \Pi_{1u}\bar \Pi_{1u}^\top + I\!\!\!\!\!}
	{\!\!\!\!\!\dfrac{m_\rho}{m_\rho^2 + \omega^2} \,\bar \Pi_{1u} \Pi_{2u}^\top\!\!\!}
	{\!\!\!*\!\!\!\!\!}
	{\!\!\!\!\!\dfrac{m_\rho/\mu}{m_\rho^2 + \omega^2} \,\Pi_{2u}\Pi_{2u}^\top + \Gamma_\rho(\omega)\!\!\!}
	\!\succeq\!
	0.
\end{align}
where
\begin{multline*}
\Gamma_\rho(\omega):=(\omega^2I + (\mu_\rho I+\gamma Q)^2)^{-1}\\
(\omega^2I+ (\mu_\rho I+\gamma Q)(\gamma Q -\rho I)).
\end{multline*}
When $\rho =0$,
\begin{multline*}
\Gamma_0(\omega):=(\omega^2I + (\mu I+\gamma Q)^2)^{-1}
(\omega^2I+ \mu \gamma Q+\gamma ^2Q^2).
\end{multline*}
Since $\mathrm{null}(\Gamma_0(\omega))\subseteq\mathrm{null}(Q)$ for all $\omega\in\R$, and $\mathrm{null}(\Pi_{2u}\Pi_{2u}^\top)\cap\mathrm{null}(Q)=\{0\}$, the bottom-right block of~\eqref{FDI5} is positive definite. Using the same Schur complement argument as in the proof of Theorem~\ref{thm.time.scale.separation}, we conclude that~\eqref{FDI5} holds \emph{strictly} for $\rho=0$. By continuity, there exist $\rho^\star>0$ such that~\eqref{FDI5} holds strictly for all $\rho<\rho^\star$. With a similar continuity argument, there exist $\epsilon^\star>0$ such that, for all $\epsilon<\epsilon^\star$,~\eqref{FDI.epsilon} is satisfied with $0\le\rho<\rho^\star$, which completes the proof.
\end{proof}

\section{Power systems example}\label{sec.example}

In this section, we provide illustrative numerical results based on a power-systems case study that is describe next.

\subsection{Online OPF for constraint-aware grid rebalancing}\label{sec.power.systems.alg}

As an example for the framework described above, we consider the linearized swing equation
\begin{align}\label{eq.swing}
\begin{split}
\begin{bmatrix}
\dot \theta\\
\dot \omega
\end{bmatrix} & = 
\overbrace{\begin{bmatrix}
0 & I\\
-M^{-1} Y & -M^{-1} D 
\end{bmatrix}}^{\bar A}
\begin{bmatrix}
 \theta\\
 \omega
\end{bmatrix} + 
\overbrace{\begin{bmatrix}
0 \\
B_1 
\end{bmatrix}}^{\bar B_u} u 
+ 
\overbrace{\begin{bmatrix}
0 \\
B_2 
\end{bmatrix}}^{\bar B_w} w \\
\underbrace{\begin{bmatrix}
 p_l\\
 \bar\omega
\end{bmatrix} }_{y}
 & = 
\underbrace{\begin{bmatrix}
E & 0 \\
0 & \frac{1}{n}\mathbf{1^\top}
\end{bmatrix}}_{\bar C}
\begin{bmatrix}
 \theta\\
 \omega
\end{bmatrix},
\end{split}
\end{align}
where the state $\xi := (\theta,\omega)$ is composed of the generator's (or inverters) phase angles $\theta_i$ and frequencies $\omega_i$, $Y$ is the grid's admittance matrix, $B_1\in\{0,1\}^{n\times m}$ and $B_2\in\{0,1\}^{n\times (n-m)}$ select the controllable and uncontrollable power injections respectively, $E$ is the matrix that maps voltage angles $\theta$ to line powers $p_l$ and $\bar \omega$ is the average system frequency. $M = \mathrm{diag}(m_1,...,m_n)$ is the matrix of rotational inertia coefficients and $D = \mathrm{diag}(d_1,...,d_n)$ is the damping matrix that collects the friction coefficients and primary control gains of the generators. System~\eqref{eq.swing} is marginally stable as $(\mathbf 1,0)$ is an eigenvector of $\bar A$ with zero eigenvalue. This is consistent with the rotational invariance of the power flow-equations. To eliminate the average-mode $\bar{\theta}$ from~\eqref{eq.swing} we introduce the following coordinate transformation~\cite{wu2016input},
	\be
	\xi
	\; = \,
	\tbo{\theta}{\omega}
	\, = \,
	\underbrace{\tbt{U}{0}{0}{I}}_{T}
	x
	~ + \;
	\tbo{\one}{0}
	\bar{\theta},
	\label{eq.x-xi}
	\ee
where the columns of the matrix $U \in \mathbb{R}^{N \times (N-1)}$ form an orthonormal basis that is orthogonal to $\textup{span} \, (\one)$. For example, these columns can be obtained from the $(N-1)$ eigenvectors of the matrix $Y$ that correspond to the non-zero eigenvalues. In the new set of coordinates,
	$
	x(t) = T^\top \xi(t) \in \bbR^{n-1},
	$
the closed-loop system takes the form
    \be
    \ba{rcl}
    \dot{x}
    & \!\! = \!\! &
     {A} \, x
    \; + \;
    {B}_u \, u
     \; + \;
    {B}_w \, w \\[0.15cm]
    y
    & \!\! = \!\! &
    {C},
    \ea
    \label{eq.ss2cl}
    \ee
where
${A} := T^\top \bar{A} \, T,
	 B_u := T^\top \bar{B}_u,
	 B_w:= T^\top \bar {B}_w$, and 
	 $C:= \bar C \, T.$
System~\eqref{eq.ss2cl} is Hurwitz stable and represents the dynamics of the rotor angles w.r.t  the average angle and the unaltered rotor frequencies. The outputs are unaltered by the coordinate transformation: $y$ remains expression of the line powers and the average grid frequency. Given $M\succ 0$, $\underline p\le\bar p$,  consider the following time-varying optimization problem
\begin{align}\label{prob.opf}
\begin{split}
\min_{x} 			&~u^\top M u + c^\top u\\
\mathrm{s.t.}  	
		  	&~\underline p \le p_l \le  \bar p\\
			&~\omega=0\\
			& 
\begin{bmatrix}
p_l \\
\bar \omega 
\end{bmatrix} = 
\Pi_{u} 
u + 
\Pi_{w} 
w,
\end{split}
\end{align}
with
$
\Pi_{u} 
 = -
C 
A^{-1}
B_u,
\quad
\text{and}
\quad
\Pi_{w} 
 = -
C 
A^{-1}
B_w.
$ 
Problem~\eqref{prob.opf} corresponds to a DC-OPF and the time-varying nature of the problem is given by the time-varying disturbance $w$ (uncontrollable loads and power injections from renewable sources). Note that normally the number of transmission lines is greater than the controllable injections, therefore the rank constraint in the assumptions of Theorem~\ref{thm.time.scale.separation} does not allow us to formulate the online primal dual method~\eqref{online.opt} using $g$ to impose both the line and frequency constraints as hard constraints. We therefore propose two alternatives
\begin{itemize}
\item \textbf{Soft constraints for the transmission lines}\\
We reformulate problem~\eqref{prob.opf} using soft constraints for the transmission line powers as
 \begin{align}\label{prob.opf.2}
\begin{split}
\min_{x} 			&~u^\top H u + c^\top u + h(y_1) + g(y_2)\\
\mathrm{s.t.} &
\begin{bmatrix}
y_1 \\
y_2 
\end{bmatrix} = 
\begin{bmatrix}
\Pi_{1u}\\
\Pi_{2u}
\end{bmatrix}
u + 
\begin{bmatrix}
\Pi_{1w}\\
\Pi_{2w}
\end{bmatrix}w
\end{split}
\end{align}
where $y_1 = p_l$ and $y_2=\bar \omega$. The soft-constraints function
$
h(y_1)_i = \eta \max\left\{\underline p_{\,i} - y_{1i} \,,\, y_{1i} -  \bar p_i\,,\,0\right\}^2 
$
is continuously differentiable with an $L_h$-Lipschitz gradient with $L_h=\eta$ 
and the indicator function
has Moreau envelope $M_{g\mu}(y_2) = \frac{1}{2\mu}\bar \omega^2$.$g(\cdot) = \mathbb I_{\{0\}}$ .
The online optimization algorithm~\eqref{online.opt} then becomes
\begin{align}\label{online.opt.opf.sc}
\begin{split}
\frac{1}{\epsilon} \dot u &= -Hu - c - \eta\Pi_{1u}^\top s_{\underline p,\bar p}(p_l) -\frac{1}{\mu}\Pi_{2u}^\top (\bar \omega+\mu\lambda) \\
\frac{1}{\epsilon} \dot \lambda &= \bar \omega, 
\end{split}
\end{align}
where $\eta\ge 0$ is a control parameter, $\epsilon$ dictates the control speed, and $s_{\underline p,\bar p}(y_1)$ is the \emph{soft thresholding} operator defined as
\[
s_{\underline p,\bar p\,i}(v) = \left\{
\ba{lc}
v_i - \underline p_i ~&~ v_i \le \underline p_i \\ 
0 ~&~ \underline p_i< v_i  <\bar p_i\\
v_i-\bar p_i ~&~ v_i \ge \bar p_i.
\ea\right.
\]

\item \textbf{Approximate online optimization}\\
A second possibility is to keep problem~\eqref{prob.opf} as it is 
 \begin{align}\label{prob.opf.3}
\begin{split}
\min_{x} 			&~u^\top H u + c^\top u  + g(y)\\
\mathrm{s.t.} &\begin{bmatrix}
y_1 \\
y_2 
\end{bmatrix} = 
\begin{bmatrix}
\Pi_{1u}\\
\Pi_{2u}
\end{bmatrix}
u + 
\begin{bmatrix}
\Pi_{1w}\\
\Pi_{2w}
\end{bmatrix}w
\end{split}
\end{align}
where $y = (y_1,y_2,) =  (p_l,\bar \omega)$  
and $g(\cdot)$ is the indicator function
$g(\cdot) = \mathbb I_{[\underline p_1,\bar p_1]\times...\times[\underline p_m,\bar p_m]\times \{0\}}$
has Moreau envelope 
\[
M_{g\mu}(y) =\frac{1}{2\mu}\sum_{i=1}^m \max\left\{\underline p_{\,i} - p_{li} \,,\, p_{li} -  \bar p_i\,,\,0\right\}^2  +  \frac{1}{2\mu}\bar \omega^2,
\] 
whose gradient is given by
\[
\nabla M_{g\mu}(y) =
\frac{1}{\mu}\begin{bmatrix}
s_{\underline p,\bar p}(p_l)\\
\bar \omega
\end{bmatrix}.
\]
Since $\Pi_{u}\Pi_{u} ^\top\not\succ 0$ we use the \emph{approximate} online optimization algorithm~\eqref{online.opt.gamma} with matrix $Q = \mathrm{diag}(I,0)$ in order to maintain a hard constraint for the frequency deviation. The algorithm~\eqref{online.opt.gamma} then becomes
\begin{align}\label{online.opt.opf.approx}
\begin{split}
\frac{1}{\epsilon}  \dot u = & -Hu - c - \frac{1}{\mu}\Pi_{1u}^\top s_{\underline p,\bar p}(p_l + \mu\lambda_1) \\
& -\frac{1}{\mu}\Pi_{2u}^\top (\bar \omega+\mu\lambda_2) \\
\frac{1}{\epsilon}  \dot \lambda_1 =&   \,s_{\underline p,\bar p}(p_l+\mu\lambda_1) - (\mu+\gamma) \lambda_1,\\ 
\frac{1}{\epsilon} \dot \lambda_2 =& \,\bar \omega. 
\end{split}
\end{align}
where $\gamma> 0$ is a control parameter and $\epsilon$ dictates the control speed.
\end{itemize}

\begin{remark}
The constraint $\Pi_{2u}\Pi_{2u} ^\top\succ 0$ is natural in the context of power systems. It is well known, for example, that decentralized integral control for frequency regulation is not internally stable as it can be destabilized by an infinitesimal bias~\cite{dorfler2017gather}. A choice of $C_2 = [0,I]$ makes $\Pi_{2u}=\alpha \mathbf 1 \mathbf 1 ^\top$, where $\alpha$ depends on the choice of transformation matrix $U$. An equality constraint on the frequency would produce a decentralized integral controller. This choice clearly violates the rank condition. 
\end{remark}
In the following section we test the proposed algorithms on a simple power systems test-case.

\subsection{IEEE9 Test Case}
\begin{center}
\vspace{-0.5cm}
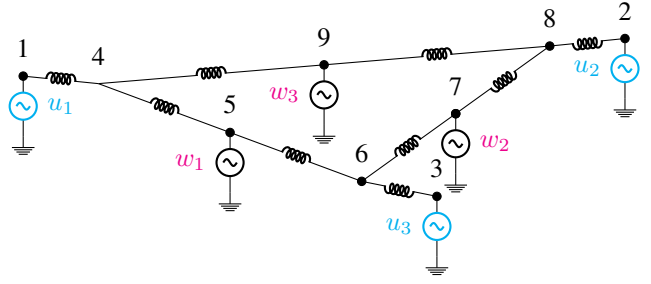
\begin{figure}
\begin{circuitikz}[american voltages]

\ctikzset{bipoles/resistor/height=0.15}
\ctikzset{bipoles/resistor/width=0.4}

\ctikzset{bipoles/generic/height=0.15}
\ctikzset{bipoles/generic/width=0.4}

\ctikzset{bipoles/length=.6cm}

\ctikzset{bipoles/thickness=2}

\node[label=1] at (0,0.1)  (I1)  {};
\node[label=2] (I2) at (8,0.6) {};
\node[label=3] (I3) at (5.5,-1.5) {};

\draw ($ (I1) + (0,-0.8) $) node [ground] (g1) {};
\draw (I1) to[sV, l=$\textcolor{ProcessBlue}{u_1}$, color=ProcessBlue] (g1);
\draw ($ (I2) + (0,-0.8) $) node [ground] (g2) {};
\draw (g2) to[sV, l=$\textcolor{ProcessBlue}{u_2}$, color=ProcessBlue] (I2);
\draw ($ (I3) + (0,-0.8) $) node [ground] (g3) {};
\draw (g3) to[sV, l=$\textcolor{ProcessBlue}{u_3}$, color=ProcessBlue] (I3);

\node[label=4] (4) at ($ (I1) + (1,-0.1) $) {};
\node[label=8] (8) at ($ (I2) + (-1,-0.1) $) {};
\node[label=6] (6) at ($ (I3) + (-1,0.2) $) {};


\node[label=5] (5) at ($(4)!0.5!(6)$) {};
\node[label=9] (9) at ($(4)!0.5!(8)$) {};
\node[label=7] (7) at ($(6)!0.5!(8)$) {};

\draw ($ (5) + (0,-0.8) $) node [ground] (g5) {};
\draw (g5) to[sV, l=$\textcolor{magenta}{w_1}$] (5);
\draw ($ (9) + (0,-0.8) $) node [ground] (g9) {};
\draw (g9) to[sV, l=$\textcolor{magenta}{w_3}$] (9);
\draw ($ (7) + (0,-0.8) $) node [ground] (g7) {};
\draw (7) to[sV, l=$\textcolor{magenta}{w_2}$] (g7);
%
%
%
%
%

\draw 		(I1) 
to[L,*-] (4);
\draw 		(I2) 
to[L,*-] (8);
\draw 		(I3) 
to[L,*-] (6);

\node (49) at ($(4)!0.5!(9)$) [label={[xshift=0cm, yshift=0.3cm]}]{};
\draw 		(4) 
to [L,-*] (9);

\node (98) at ($(8)!0.5!(9)$) [label={[xshift=0cm, yshift=0.3cm]}]{};
\draw 		(9) 
to [L,-*] (8);

\node (45) at ($(4)!0.5!(5)$) [label={[xshift=0cm, yshift=0.3cm]}]{};
\draw 		(4) 
to [L,-*] (5);

\node (56) at ($(6)!0.5!(5)$) [label={[xshift=0cm, yshift=0.3cm]}]{};
\draw 		(5) 
to [L,-*] (6);

\node (67) at ($(6)!0.5!(7)$) [label={[xshift=0cm, yshift=0.3cm]}]{};
\draw 		(6) 
to [L,-*] (7);

\node (78) at ($(8)!0.5!(7)$) [label={[xshift=0cm, yshift=0.3cm]}]{};
\draw 		(8) 
to [L,-*] (7);

%
\end{circuitikz}
\caption{Simplified IEEE9 Test case: we consider a lossless system with three controllable generators and three uncontrollable loads}
\label{fig.IEEE9}
\vspace{-0.5cm}
\end{figure}
\end{center}
In this section we study the proposed algorithms on a power-system case study. We consider a modified version of the IEEE9 test case depicted in Figure~\ref{fig.IEEE9} where we neglect the line resistances and shunt capacitances. We consider three controllable generators whose power setpoints are updated according to the online optimization algorithms outlined in Section~\eqref{sec.power.systems.alg} and three uncontrollable loads as Illustrated in Figure~\ref{fig.IEEE9}. We consider the swing dynamics~\eqref{eq.swing} with the following parameters $m_i=1$ p.u. for the controllable generators and $m_i=0.1$ for the uncontrollable loads/generators and $d_i=0.1$ p.u. for all generators and loads. We consider a ``nominal'' disturbance of $w =( -0.9,-1,-1.25)$ p.u. and we solve the OPF problem~\eqref{prob.opf} for $H = I$ and $c = \mathbf 1$. We simulate both algorithms proposed in Section~\eqref{sec.power.systems.alg} with the following events: the constraint on Line $1$ (from bus $1$ to $4$) is moved from 
$1.5$ to $0.5$ p.u. after 10 seconds and restored after 20 seconds. After $50$ seconds the $w_1$ becomes $w_1(t) = -0.9*[1+\cos(0.2\,t)]$ effectively alternating between no load and double the nominal load. 

\begin{figure}[t]
\begin{center}
\input{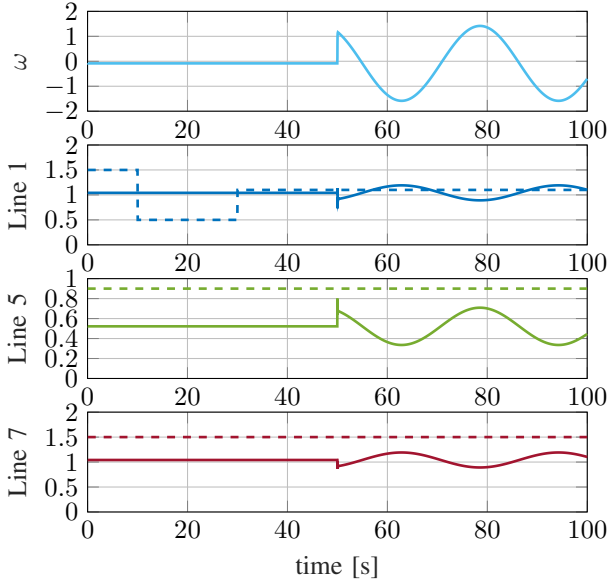}
\vspace{-.2cm}
\caption{The IEEE9 test system in the proposed scenario without control: the system does not react the the tightening of the constraint on Line 1 and cannot regulate the frequency to react to the load variation.}
\label{fig.no.control}
\vspace{-.5cm}
\end{center}
\end{figure}

\begin{figure}[t]
\begin{center}
\input{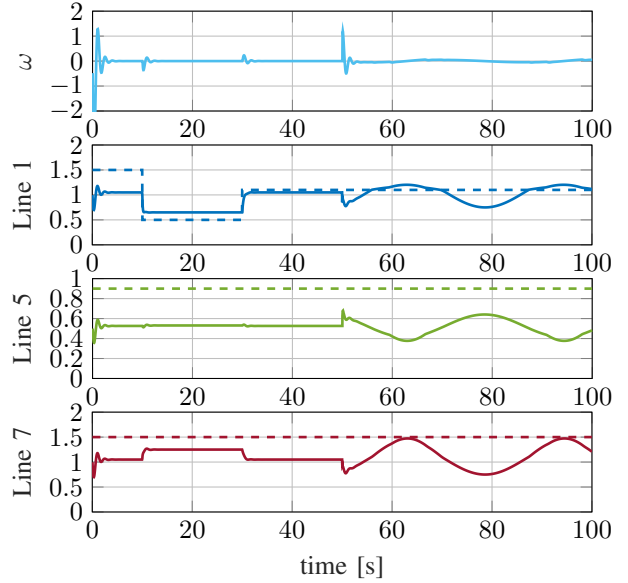}
\vspace{-.3cm}
\caption{IEEE9 test system with online optimization using soft constraints on the transmission lines~\eqref{online.opt.opf.sc} with parameters $\eta =4,\mu =4$, $\epsilon = 10^{-2}$. The system is able to react to the tightened constraint on Line 1 and keep the frequency regulated to zero. }
\label{fig.soft}
\vspace{-.5cm}
\end{center}
\end{figure}
\begin{figure}[t]
\begin{center}
\input{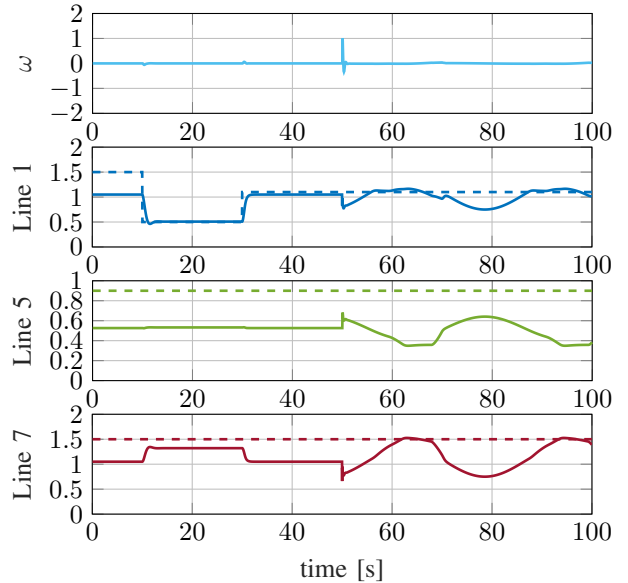}
\vspace{-.2cm}
\caption{IEEE9 test system with approximate online optimization using hard constraints and $\gamma-$regularization on the transmission lines~\eqref{online.opt.opf.approx} with parameters $\gamma = 10^-2$, $\mu =4$, $\epsilon = 10^{-2}$. The approximate algorithm has superior performance both in reacting to the tightened constraint on Line 1 and keeping the frequency regulated to zero. }
\label{fig.approx}
\vspace{-.5cm}
\end{center}
\end{figure}

In Figures~\ref{fig.no.control}-\ref{fig.approx} we illustrate the behavior of the system average frequency and line powers under different control laws. The dashed lines represent the line constraints. In Figure~\ref{fig.no.control} we observe the behavior of the system without control. Predictable, the system is unable to react to the tightening of the constraint on Line 1 and to keep the frequency and line powers within limits when the load is changing. In Figure~\ref{fig.soft} we show the system controlled with the online optimizer~\eqref{online.opt.opf.sc} with parameters $\eta =4,\mu =4$, $\epsilon = 10^{-2}$. The algorithm~\eqref{online.opt.opf.sc} uses soft constraints for the transmission lines and hard constraints for the frequency. We observe that the system keeps the frequency regulated and reacts to the constraints tightening but does not reach constraint satisfaction at steady-state. This is an artifact of the soft constraints. Finally, in Figure~\ref{fig.approx} We illustrate the behavior of the approximate online optimization algorithm~\eqref{online.opt.opf.approx} with parameters $\gamma = 10^-2$, $\mu =4$, $\epsilon = 10^{-2}$. All constraints are kept as hard constraints, but the dual update on the lagrange multipliers corresponding to the line constraints is regularized with the parameter $\gamma>0$ as illustrated in~\eqref{online.opt.opf.approx}. The algorithm is superior to the soft constrained optimization both in frequency regulation and in enforcing the constraints. The stability of both algorithms was certified with the appropriate LMI of the form~\eqref{LMI} for $\rho =10^{-3}$. Note that with the proposed framework we can include some specialized frequency control algorithms proposed in the literature~\cite{li2016connecting,dorfler2016breaking} and extend them to more general constraints. Distributed implementation the proposed algorithms is an interesting topic for future research.

\section{Conclusion and outlook}\label{sec.conclusion}

This paper studied the problem of regulating the output of a LTI dynamical system to the solution of a time varying convex optimization problem. In particular, we proposed a feedback interconnection of an augmented Lagrangian primal-dual saddle-flow algorithm with the system. We then provided LMI conditions based on an IQC characterization of the nonlinear portion of the algorithm that guarantees that the interconnection is stable and the output of the dynamical system tracks the time-varying optimizer of the optimization problem. Furthermore, we proved that, under mild conditions, the interconnection is always stable provided there is enough time-scale separation between the system (fast) and the optimizer (slow). We observed that the LMI conditions are not very conservative in establishing stability. The guaranteed rate of convergence estimate, however, tends to be extremely conservative. To this end we plan to explore more complex IQC-based approaches such as the Popov IQC or the Zames-Falbes IQC that, at least for the time-invariant case, promise to provide tighter rates of convergence. This work opens to numerous research directions: among them are the investigation of decentralized and distributed implementations of the online optimization algorithms, exploration of second order methods based on the proximal augmented Lagrangian~\cite{dhingra2017second} and the study of convergence in case of non-convex problem or nonlinear systems. 

\section*{Acknowledgements}
We kindly acknowledge Neil K. Dhingra and Prof. Mihailo R. Jovanovi\'c for fruitful discussions on proximal augmented Lagrangian methods.

\bibliographystyle{unsrt}
\bibliography{bib_file,biblio_part2}

\end{document}